\documentclass[final,5p,twocolumn,times]{elsarticle}

\usepackage{graphicx}          % Include this line if your 
\usepackage{mathptmx} % assumes new font selection scheme installed
\usepackage{times} % assumes new font selection scheme installed
\usepackage{amsmath} % assumes amsmath package installed
\usepackage{amssymb}  % assumes amsmath package installed

\usepackage{balance}
\usepackage{natbib}
\usepackage{color}

\newtheorem{remark}{\textbf{Remark}}
\newtheorem{theorem}{\textbf{Theorem}}

\begin{document}

\begin{frontmatter}


\title{\LARGE Dead-beat model predictive control for discrete-time linear systems}  % \tnoteref{Funds}} 
												% Title, preferably not more 

%\thanks{%This paper was not presented at any IFAC meeting. 
%Corresponding author B.~Zhu. Tel. +86-13466314059. 
%Fax +86-10-82339382.}

\tnotetext[Funds]{%This paper was not presented at any IFAC conferences.
This work was supported by National Natural Science Foundation of China under grants 62073015.
%Corresponding author B.~Zhu. Tel. +86 13466314059.
%Fax +86 10 82339382.
}

\author[buaa]{Bing Zhu}\ead{zhubing@buaa.edu.cn}    % Add the 
%\author[Rome]{Julius Caesar}\ead{julius@caesar.ir},               % e-mail address 
%\author[Baiae]{Publius Maro Vergilius}\ead{vergilius@culture.ir}  % (ead) as shown

\address[buaa]{The Seventh Research Division,
	School of Automatic Science and Electrical Engineering,\\
	Beihang University, Beijing 100191, P.R.~China}  % Please supply                                              
%\address[Rome]{Senate House, Rome}             % full addresses
%\address[Baiae]{The White House, Baiae}        % here.

\begin{keyword}                           % Five to ten keywords,  
Model predictive control; dead-beat control; linear systems; constrained control              
% chosen from the IFAC 
\end{keyword}                             % keyword list or with the 
                                          % help of the Automatica 
                                          % keyword wizard

\begin{abstract}                          % Abstract of not more than 200 words.
In this paper, 
model predictive control (MPC) strategies are proposed for dead-beat control of linear systems with and without state and control constraints.
In unconstrained MPC, deadbeat performance can be guaranteed by setting the control horizon to the system dimension, and adding an terminal equality constraint.
It is proved that the unconstrained deadbeat MPC is equivalent to linear deadbeat control.
The proposed constrained deadbeat MPC is designed by setting the control horizon equal to the system dimension and penalizing only the terminal cost.
The recursive feasibility and deadbeat performance are proved theoretically.
\end{abstract}

\end{frontmatter}

\section{Introduction}

For discrete-time systems, deadbeat control implies that the closed-loop state is identically zero after finite time steps.
The essence of deadbeat control for SISO discrete-time linear systems was discovered in \cite{Kalman1960On},
where existence of linear deadbeat control was proved to be equivalent to controllability.
A well known design strategy for linear deadbeat control is to assign all eigenvalues to the origin.
More fundamental and systematic results of linear deadbeat control can be found in \cite{o1981discrete, kucera1984deadbeat} and references therein.
%From engineering viewpoint, sometimes only the output is required to have deadbeat pefermance,
%and typical results on output deadbeat control can be found in \cite{bastin1999output, marrari1989output} and references therein.

Deadbeat control can be designed in optimal control framework,
where the core idea is calculating the deadbeat gain from quadratic cost function with proper weighting matrices and boundary conditions.
In \cite{Leden1976Dead}, the deadbeat control was solved from singular Riccati equation, 
where %only terminal cost was penalized in the cost function.
recursive strategy was applied to solve the singular Riccati equation.
Analytical solution to singular Riccati equation was provided in \cite{Lewis1981A}, 
where  %Moore-Penrose inverse is applied, and 
deadbeat control was considered as a special case.
Another strategy for optimal deadbeat control is to only penalize the state cost without penalizing the control cost \citep{emami1982deadbeat}.
%This strategy can be further extended to deal with multiple transmission zeros \citep{spurgeon1991output}. 
In \cite{sugimoto1993direct}, a direct solution is proposed for singular Riccati equation to calculated the deadbeat gain.

%In its early years, 
The deadbeat control theory has been developed maturely.
Recent results on deadbeat control usually focus on engineering applications.
Deadbeat control is applicable especially to systems requiring high accuracy and short transient process \citep{meng2022optimized}.
Deadbeat strategy is also useful in observer design, 
where exact zero estimation error is guaranteed, 
such that overall performances with estimated states can be improved 
\citep{Meng2022feedback, masti2021machine, shang2021distributionally}.
For nonlinear systems, deadbeat control is a relatively open problem, 
and some results are provided in \cite{haddad2020finite, tuna2012state} and references therein.

In some cases, the systems are subject to hard constraints,
and bounded deadbeat control should be considered.
Bounded deadbeat control problem was firstly investigated geometrically in \cite{Wing1963THE},
where the plant eigenvalues were required to be inside the closed unite circle to guarantee the global admissibility.
Model predictive control (MPC) can be introduced to guarantee bounded deadbeat control for some specific systems,
e.g., permanent-magnet synchronous motor \citep{kang2019symmetrical}, interconnected heterogenerous power converters \citep{kreiss2021optimal}, and
electrical vehicle charging systems \citep{wang2019improved}, etc.
Currently, however, fundamental and systematic results on deadbeat MPC theory are relatively limited.

In this paper, we are proposing fundamental solutions to the deadbeat control problem in the MPC framework,
such that the deadbeat performance can be achieved by bounded control with improved transient process.
\textit{Main contributions} of this paper include:
1) deadbeat controls are designed and proved in a systematic MPC framework;
2) implicit and explicit solutions are provided for unconstrained linear deadbeat MPC;
and 3) constrained deadbeat MPC is designed for linear systems subject to state and control constrains.

%==========================================================================
\section{Problem formulation}

%\subsection{Problem statement}

In this paper, we consider the
single-input-single-output (SISO) discrete-time linear system:
\begin{align}\label{linear syst}
 x(k+1) = Ax(k)+Bu(k)
\end{align}
where $x\in\mathbb{R}^n$ and $u\in\mathbb{R}$  are the system state and control input, respectively;
%the state matrix $A$ and the input matrix $B$ have proper dimensions,
and $(A,B)$ is controllable.

The system is subject to state and control constraints
\begin{align}\label{constraints}
	x\in \mathcal{X},~~u\in \mathcal{U},
\end{align}
where $\mathcal{X}$ and $\mathcal{U}$ are convex, and contain the origin.

The objective is to design $u(k)$ in MPC frameworks such that, 
for some $T>0$, it holds that $x(k)=0$ for all $k\geq T$.

%===============================================================

\section{Unconstrained dead-beat MPC}

In this section, no constraints are exerted on states or control inputs,
i.e., constraints \eqref{constraints} is not considered in this section.
%Usually, the unconstrained MPC is equivalent to linear feedback control.
%Theoretically, by using linear feedback control,
%the closed-loop system is capable of reaching the origin within steps no more than its dimension.
%
In the framework of unconstrained MPC, deadbeat control can be implemented through two strategies,
namely deadbeat MPC with terminal equality constraint and deadbeat MPC with only terminal cost.
Both strategies requires the control horizon be equal to state dimensions, namely $N=n$.

%------------------------------------------------------------------------
\subsection{Uncosntrained deadbeat MPC with terminal equality constraint}
\label{sec unconstrained}

In this deadbeat strategy, we add an extra terminal equality constraint to the original unconstrained optimization to guarantee the finite-time convergence.

Define the state sequence and control sequence at time $k$ by
\begin{align*}
	X(k) =& [x^T(1|k), \cdots, x^T(N|k)]^T,\\
	U(k) =& [u^T(0|k), \cdots, u^T(N-1|k)]^T,
\end{align*}
where $(i|k)$ denotes the prediction $i$ steps ahead from time $k$;
and $N$ is the control horizon.
Throughout the paper, the prediction horizon is set equal to the control horizon.

In this deadbeat MPC, 
two key settings include: 1) the control horizon is set equivalent to state dimension, i.e., $N=n$;
and 2) a terminal equality constraint $x(n|k) = 0$ is added to the original unconstrained optimization.

The optimization is formulated by
\begin{align}
	&U^\ast(k) = \mathrm{arg}\min_{U(k)} \sum_{i=1}^n\|x(i|k)\|^2_Q+\|u(i-1|k)\|^2_R,  \label{unconstrained opt1}\\
	\mbox{s.t.}    ~~ &x(i+1|k) = Ax(i|k)+Bu(i|k),~~i=0,\cdots,n-1, \\
	&x(n|k) = 0, \label{terminal equality constraint}
\end{align}
where $n$ is the dimension of $x$.
The deadbeat MPC is implemented in receding horizon scheme:
\begin{align}\label{receding horizon scheme}
 u(k) = [1~0~\cdots~0]U^\ast(k).
\end{align}

%\begin{remark}
%	Although the system is supposed to be unconstrained,
%	the optimization in deadbeat MPC is actually constrained.
%\end{remark}

\begin{theorem}\label{deadbeat thm}
	Consider the linear system \eqref{linear syst}, and suppose that no constraints are exerted.
	Suppose also that $(A,B)$ is completely controllable.
	The MPC is calculated from the optimization \eqref{unconstrained opt1}--\eqref{terminal equality constraint},
	and implemented by \eqref{receding horizon scheme}.
	Then, 
	\begin{enumerate}
	\item[1)]
	the optimization \eqref{unconstrained opt1}--\eqref{terminal equality constraint} is feasible globally;
	\item[2)]
	the closed-loop state satisfies $x(k)=0$ for all $k>n$.
	\end{enumerate}
\end{theorem}

\emph{\textbf{Proof}}:
The behavior of \eqref{linear syst} can be predicted by
\begin{align}
	x(1|k) =& Ax(k) +Bu(0|k),\\
	x(2|k) =& A^2x(k) +ABu(0|k)+Bu(1|k),\\
	\vdots & \nonumber\\
	x(n|k) = & A^{n}x(k) + [A^{n-1}B~A^{n-2}B~\cdots~B]U(k), \label{predictive state sequence}
\end{align}
where, since $(A,B)$ is controllable, the square matrix
\begin{align*}
	S\triangleq [A^{n-1}B~A^{n-2}B~\cdots~B]
\end{align*}
is invertible.
It indicates that a unique control sequence $U^\ast(k)$ always exists,
such that \eqref{terminal equality constraint} is satisfied.
Consequently, the optimization is always feasible regardless of initial states. This completes the proof of 1).

At time $k=0$, the optimal control sequence is unique, and the corresponding optimal state sequence is also unique:
\begin{align*}
	U^\ast(0) =& [u^\ast(0|0), \cdots, u^\ast(n-1|0)]^T,\\
	X^\ast(0) =& [x^{\ast T}(1|0), \cdots, x^{\ast T}(n-1|0), 0]^T,
\end{align*}
where $u(0) = u^\ast(0|0)$ is implemented, and $x(1) = x^{\ast}(1|0)$ can be obtained.

At time $k=1$, the optimal control sequence and the optimal state sequence are still unique, and they have to be
\begin{align*}
	U^\ast(1) =& [u^\ast(1|0), \cdots, u^\ast(n-1|0), 0]^T,\\
	X^\ast(1) =& [x^{\ast T}(2|0), \cdots, x^{\ast T}(n-1|0), 0, 0]^T,
\end{align*}
where $u(1) = u^\ast(0|1)=u^\ast(1|0)$ is implemented, and $x(2) = x^{\ast}(1|1) = x^\ast(2|0)$ can be obtained.

At time $k=n-1$, the unique control sequence and the optimal state are unique, and they have to be
\begin{align*}
	U^\ast(n-1) =& [u^\ast(n-1|0) ,0, \cdots, 0]^T,\\
	X^\ast(n-1) =& [0,\cdots, 0]^T,
\end{align*}
where $u(n-1) = u^\ast(0|n-1)=u^\ast(n-1|0)$ is implemented, and $x(n) = 0$.
It implies that, $u(k) = 0$ and $x(k) = 0$ for all $k>n$. This completes the proof of 2).
\hfill $\square$

It follows from the terminal equality constraint \eqref{terminal equality constraint} 
and the predictive state sequence \eqref{predictive state sequence} that
the control sequence can be solved uniquely and explicitly by
\begin{align}\label{explicit control sequence}
	U(k) = - S^{-1}A^{n}x(k),
\end{align}
and the MPC can be calculated explicitly by
\begin{align}
	u(k) = -[1~0~\cdots~0]S^{-1}A^{n}x(k) = -K_{db}x(k),
\end{align}
where the deadbeat control gain is calculated by
\begin{equation}\label{explicit gain}
	K_{db} = [1~0~\cdots~0]S^{-1}A^{n} = S^T_nA^n,
\end{equation}
and $S_n^T$ is the first row of $S^{-1}$.

\begin{theorem}\label{thm exp}
	If the control gain is calculated by \eqref{explicit gain}, then all eigenvalues of $A-BK_{db}$ are zero,
	and $x(k)=0$ for all $k>n$.
\end{theorem}

\emph{\textbf{Proof}}:
For matrix $A$, there always exists invertible transformation matrix $T$, such that
\begin{align*}
	T AT^{-1} = \begin{bmatrix}
		0_{(n-1)\times 1} & I_{(n-1)\times(n-1)}\\ -a_0 &-a_1 ~~\cdots ~~-a_{n-1}
	\end{bmatrix},
~~TB = \begin{bmatrix}
	0_{(n-1)\times 1}\\ 1
\end{bmatrix}.
\end{align*}
According to linear control theory \citep{Ogata2010Modern}, 
the tranformation matrix $T$ is calculated by
\begin{align*}
	T = \begin{bmatrix}
		S_n & A^TS_n &\dots & (A^{n-1})^TS_n
	\end{bmatrix}^T.
\end{align*}
%Here, it is desirable to prove
%\begin{align*}
%	TBK_{db}T^{-1} = \begin{bmatrix}
%		0 &0 &0 &\cdots &0\\  \vdots & &0 & &\vdots\\
%		& & &\ddots &0\\ 0 &0 &\cdots &0 &0\\ -a_0 &-a_1 &\cdots & &-a_{n-1}
%	\end{bmatrix},
%\end{align*}
%such that 
%\begin{align*}
%	T(A-BK_{db})T^{-1} = \begin{bmatrix}
%		0 &I_{n-1}\\ 0 &0
%	\end{bmatrix},
%\end{align*}
%and eigenvalues of $A-BK_{db}$ are all zero.
It then follows that
\begin{align*}
	TBK_{db}= \begin{bmatrix}
		0_{(n-1)\times 1}\\ 1
	\end{bmatrix}[1~0~\cdots~0]S^{-1}A^n
	=\begin{bmatrix}
		0_{(n-1)\times n}\\ S_n^TA^n
	\end{bmatrix},
\end{align*}
where $S_n^T$ is the first row of $S^{-1}$.
In another aspect,
\begin{align*}
	\begin{bmatrix}
		0_{(n-1)\times n} \\  -a_0 ~\cdots ~-a_{n-1}
	\end{bmatrix}T =& 
\begin{bmatrix}
	0_{(n-1)\times n} \\  -S_n^T\sum_{i=0}^{n-1}a_iA^i
\end{bmatrix} = \begin{bmatrix}
0_{(n-1)\times n}\\ S_n^TA^n
\end{bmatrix}    \nonumber\\
=& TBK_{db}, 
\end{align*}
where Cayley-Hamilton theorem is applied.
Consequently,
\begin{align*}
	TBK_{bd}T^{-1} =& \begin{bmatrix}
		0_{(n-1)\times n} \\  -a_0 ~\cdots ~-a_{n-1}
	\end{bmatrix},\\
%\end{align*}
%such that 
%\begin{align*}
	T(A-BK_{db})T^{-1} =& \begin{bmatrix}
			0_{(n-1)\times 1} &I_{(n-1)\times(n-1)}\\ 0 &0_{1\times(n-1)}
		\end{bmatrix},
\end{align*}
indicating eigenvalues of $A-BK_{db}$ are all zero.
As a result, it holds that $x(k)=0$ for all $k>n$.
\hfill $\square$

%\begin{remark}
%	The deadbeat MPC with terminal equality constraint cannot be extended to constrained cases.
%	The reason is that, in case of (hard) state or control constraints, 
% 	the unique solution $U^\ast(0)$ is not guaranteed.
%\end{remark}

%---------------------------------------------------------------

\subsection{Unconstrained dead-beat MPC with terminal cost}

This deadbeat strategy is cited from \cite{Lewis1981A},
where the deadbeat control is solved from inverse optimization.
%where it is a special case of the inverse solution to unconstrained optimization.
We explain the deadbeat control in MPC framework,
and extended it to the constrained case in the next section.

The optimization only penalizes the terminal cost without penalizing other stage costs.
In this regard, the cost function is constructed by
$
	J(k) = x^T(n|k)Px(n|k)
$,
%where $P$ is the positive definite weighting matrix, 
and the control horizon is set to the state dimension.
%No other extra constraints are required here.

The optimization is established by
\begin{align}\label{unconstrained opt terminal cost}
	&U^\ast(k) = \mathrm{arg}\min_{U(k)} x^T(n|k)Px(n|k),\\
	\mbox{s.t.}~~&x(i+1|k) = Ax(i|k)+Bu(i|k),~~i=0,\cdots,n-1,
\end{align}
and the deadbeat MPC is implemented by \eqref{receding horizon scheme}.

\begin{theorem}
	Consider the linear system \eqref{linear syst}, and suppose that no constraints are exerted.
	Suppose also that $(A,B)$ is completely controllable.
	The MPC is calculated from optimization \eqref{unconstrained opt terminal cost}
	and implemented by \eqref{receding horizon scheme}.
	Then, 
%	\begin{enumerate}
%		\item[1)]
%		the optimization \eqref{unconstrained opt terminal cost} is feasible globally;
%		\item[2)]
		the closed-loop state satisfies $x(k)=0$ for all $k>n$.
%	\end{enumerate}
\end{theorem}

\textit{\textbf{Proof}}:
%Since $(A,B)$ is completely controllable,  
%there always exists a feasible and unique control sequence $U^\ast$ driving the state from any $x(0)$ to $x(n|k)$ in $n$ steps. This proves 1).
Without any hard constraints,
the optimal terminal state minimizing $x^T(n|k)Px(n|k)$ is $x(n|k)=0$.
Then the proof follows from that of Theorem \ref{deadbeat thm}.
\hfill $\square$

The optimization can also be solved analytically by
\begin{align}\label{partial derivative}
	\left.\frac{\partial \left(x^T(n|k)Px(n|k)\right)}{\partial U} \right|_{U=U^\ast}= 0,
\end{align}
where $x(n|k)$ is calculated by \eqref{predictive state sequence}.
Solving \eqref{partial derivative} yields
\begin{align}\label{solution of partial derivative}
	2(A^nx(k)+SU^\ast(k))^TPS=0,
\end{align}
where $P$ and $S$ are invertible. 
The explicit solution to \eqref{solution of partial derivative} is exactly the same as \eqref{explicit control sequence}--\eqref{explicit gain}.
%and the deadbeat result follows from Theorem \ref{thm exp}.

%\begin{remark}
%	It can be seen that the previous two deadbeat unconstrained MPC are equivalent.
%\end{remark}

%===============================================

\section{Constrained dead-beat MPC}

The unconstrained deadbeat MPC with terminal cost can be extended to the constrained case.
Key settings include:
1) the control horizon is set equal to the dimension of system state, i.e., $N = n$;
2) only the terminal cost is penalized, i.e., $J(k)=x^T(n|k)Px(n|k)$; and
3) the postive definite weighting matrix $P$ is solved from the Lyapunov equation
\begin{align}\label{Lyapunov equation}
	A_K^TPA_K - P = -Q - K^TRK,
\end{align}
where $Q$ and $R$ are positive definite matrices;
the control gain $K$ is selected such that all eigenvalues of $A_K=A-BK$ are inside the unit circle.

The optimization is constructed by
\begin{align}
	&U^\ast(k) = \mathrm{arg}\min_{U(k)}x^T(n|k)Px(n|k),  \label{constrained opt}\\
%\end{align}
%\begin{align}
\mbox{s.t.}~~&x(i+1|k) = Ax(i|k)+Bu(i|k),~~i=0,\cdots,n-1,\\
	&x(i|k)\in\mathcal{X},~~u(i-1|k)\in\mathcal{U},~~i=1,\cdots,n,   \label{stage constraints}\\
	&x(n|k)\in \mathcal{X}_f,		\label{terminal constraints}
\end{align}
where $\mathcal{X}_f$ denotes the terminal constraint satisfying
\begin{align}
	&\mathcal{X}_f\subset \mathcal{X},\\
	&x(k)\in \mathcal{X}_f~~\Rightarrow~~ -K_{db}x(k)\in \mathcal{U},   \label{inv control constraint}\\
	&x(k)\in \mathcal{X}_f~~\Rightarrow~~ (A-BK_{db})x(k) \in \mathcal{X}. \label{inv state constraint}
\end{align}
The MPC is impemented in receding horizon scheme \eqref{receding horizon scheme}.

\begin{remark}
	The condition \eqref{inv state constraint} indicates that $\mathcal{X}_f$ is an invariant set of the deadbeat closed-loop system.
	Here we do not have to calculate $K_{db}$ in advance to obtain the largest invariant set;
	instead, a sufficiently small set $\mathcal{X}_f$ can be chosen such that it is a subset of the largest invariant set.
\end{remark}

\begin{theorem}
	Consider system \eqref{linear syst} subject to constraints \eqref{constraints}.
	The MPC is calculated by the constrained optimization \eqref{constrained opt}--\eqref{terminal constraints}, 
	and implemented by \eqref{receding horizon scheme},
	where the terminal weighting matrix is calculated from \eqref{Lyapunov equation}.
	Suppose that the constrained optimization \eqref{constrained opt} is feasible initially.
	Then,
	\begin{enumerate}
		\item[1)] constrained optimization \eqref{constrained opt} is feasible recursively;
		\item[2)] the closed-loop system is asymptotically stable;
		\item[3)] there exists a finite time $T>0$, 
		such that the state of the closed-loop system satisfies $x(k)=0$ for all $k>T$.
	\end{enumerate}
\end{theorem}

\textbf{\textit{Proof}}:
1)
%Proof of 1) is straightford in a typical MPC way by repeating the tail of the prvious step.
Suppose, at time $k$, the constrained optimization is feasible, i.e.
$U^\ast(k)$ exists such that \eqref{constraints} and \eqref{terminal constraints} are satisfied.

At time $k+1$, at least one feasible control sequence exists:
\begin{align*}
	&u(i|k+1) = u^\ast(i+1|k), ~~i = 0,\cdots, n-1,\\
	&u(n-1|k+1) = -K_{db}x^\ast(n|k),
\end{align*}
where $u^\ast(i-1|k)$ denotes the optimal control sequence at time $k$,
and $x^\ast(i|k)$ denotes the corresponding optimal state sequence at time $k$.

It is clear that $u(i|k+1),~i = 0,\cdots, n-1$ satisfy the constrol constraint,
and $u(n-1|k+1) = -K_{db}x^\ast(n|k)$ satisfies the control constraint due to \eqref{inv control constraint}.
It is also clear that $x(i|k+1),~i = 1,\cdots, n$ satisfy the state constraint,
and $x(n|k+1)=(A-BK_{db})x(n-1|k+1) = (A-BK_{db})x^\ast(n|k)$ satisfies the state constraint due to \eqref{inv state constraint}.

Consequently, at time $k+1$, the constrained optimization is feasible, 
provided that it is feasible at time $k$.

2) Take $J^\ast(k)= x^{\ast T}(n|k)Px^\ast(n|k)$ as the Lyapunov candidate for $x(n|k)$.
It follows that
\begin{align*}
	&J^\ast(k+1)-J^\ast(k) \leq J(k+1) - J^\ast(k)\\
	= & x^{T}(n|k+1)Px(n|k+1) -x^{\ast T}(n|k)Px^\ast(n|k)\\
	= & x^{\ast T}(n|k)A_K^TPA_Kx^\ast(n|k)-x^{\ast T}(n|k)Px^\ast(n|k)\\
	=& x^{\ast T}(n|k)\left(A_K^TPA_K-P\right)x^\ast(n|k)\\
	=& -\|x^\ast(n|k)\|^2_Q-\|K_{db}x^\ast(n|k)\|^2_R,
\end{align*}
where $A_K = A-BK$.
It indicates that $x^\ast(n|k)\rightarrow 0$ as $k\rightarrow \infty$.

The terminal state $x^\ast(n|k)$ at each time $k$ satisfies
\begin{align}\label{terminal state}
	x^\ast(n|k) = A^nx(k)+SU^\ast(k),
\end{align}
where $S=[A^{n-1}B, A^{n-2}B,\cdots,B]$ is invertible.
It follows from \eqref{terminal state} that, for each $x^\ast(n|k)$,
the constrained control sequence is unique:
%\begin{align*}
	$U^\ast(k) = S^{-1}\left(x^\ast(n|k)-A^nx(k)\right)$,
%\end{align*}
and the control is implemented by
\begin{align}
	u^\ast(k) = &[1, 0, \cdots, 0]S^{-1}\left(x^\ast(n|k)-A^nx(k)\right),\\
	=& -K_{db}x(k)+S_n^Tx^\ast(n|k),    \label{actual control}
\end{align}
where $K_{db}$ is the unconstrained deadbeat gain obtained in Section \ref{sec unconstrained},
and $S_n^T = [1, 0, \cdots, 0]S^{-1}$ is the first row of $S^{-1}$.

Substituting the control \eqref{actual control} into the original system yields
\begin{align}\label{perturbed}
	x(k+1) = (A-BK_{db})x(k)+BS_n^Tx^\ast(n|k),
\end{align}
which is the closed-loop deadbeat control system perturbed by the vanishing term $BS_n^Tx^\ast(n|k)$.

The solution to \eqref{perturbed} is calculated by
\begin{align*}
	x(k) = A^k_{db}x(0)+[A_{db}^{k-1}B, \cdots, B]\begin{bmatrix}
		S_n^Tx^\ast(n|0)\\ \vdots\\ S_n^Tx^\ast(n|k)
	\end{bmatrix},
\end{align*}
where $A_{db}=A-BK_{db}$.

Whenever $k\geq n$, it holds that $A_{db}^k=0$, and
\begin{align}\label{solution}
	x(k) = [A_{db}^{n-1}B, \cdots, B]\begin{bmatrix}
		S_n^Tx^\ast(n|k-n+1)\\ \vdots\\ S_n^Tx^\ast(n|k)
	\end{bmatrix}.
\end{align}
Since $x^\ast(n|k)\rightarrow 0$ as $k\rightarrow\infty$,
it follows that $x(k)\rightarrow 0$ as $k\rightarrow\infty$,
i.e., the closed-loop system is asymptotically stable.

3) Based on 1) and 2), $x(k)$ converges to the origin asymptotically,
indicating that it enters $\mathcal{X}_f$ in finite time.

Inside the invariant set $\mathcal{X}_f$, the state and control constraints are actually inactive,
and the optimization is actually unconstrained.
Consequently, the behabior of $x(k)$ follows from that of the unconstrained deadbeat MPC,
and it reaches the origin in finite time and maintains zero thereafter.
\hfill $\square$

%\begin{remark}
%	It is not required that the control gain $K$ be chosen equal to $K_{db}$.
%\end{remark}

%===============================================================
\section{Numerical example}

The plant to be controlled is given by
\begin{align*}
	A = \begin{bmatrix}
		1.1 &2 &0\\ 0 &0.95 &1\\ 0 &0 &1.2
	\end{bmatrix},~~
	B = \begin{bmatrix}
		0\\ 0.079\\ 0.1
	\end{bmatrix},
\end{align*}
where its dimension is $n=3$.
%It is straightforward to verify that $(A,B)$ is completely controllable.

%-------------------------------------------------------------------
\subsection{Unconstrained deadbeat MPC}
\label{sec sim unconstrained eq}

In this section, the unconstrainned deadbeat MPC can be calculated numerically by \eqref{unconstrained opt1} and \eqref{receding horizon scheme},
where the positive weighting matrices $Q= I_{3\times 3}$ and $R=0.1$ are set;
the control horizon is set to $N=n=3$; and the terminal equality constraint is assigned to
$x(3|k) = 0$.

\begin{figure}
%\begin{center}
\includegraphics[scale=0.35]{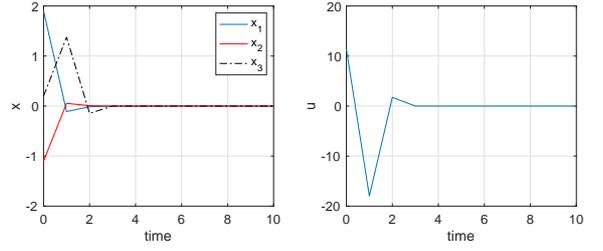}    
\caption{States and control of the closed-loop system with implicit unconstrained deadbeat MPC using the terminal equality constraint:
they reach the origin in $3$ steps.}  
\label{fig: unconstrained terminal eq imp}                             
%\end{center}                                
\end{figure}

The simulation results are displayed in Fig.~\ref{fig: unconstrained terminal eq imp}, 
where it can be seen that system states and the control input reach the origin in $3$ steps and are maintained thereafter.
To achieve deadbeat performance, the magnitude of control might be relatively large in the transient process.

\begin{figure}
\begin{center}
\includegraphics[scale=0.35]{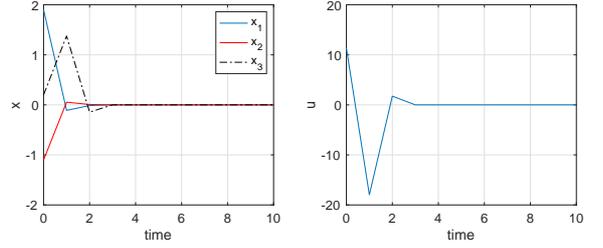}    
\caption{States and control of the closed-loop system with explicit unconstrained deadbeat MPC using the terminal equality constraint:
performances are the same as those in Fig.~\ref{fig: unconstrained terminal eq imp} .} 
\label{fig: unconstrained terminal eq exp}                             
\end{center}                                
\end{figure}

The explicit deadbeat MPC gain is calculated from \eqref{explicit gain} by
\begin{align*}
	K_{db} =& [1~0~0][A^2B~AB~B]^{-1}A^3 
	=[7.2258  ~25.1192  ~12.6558],
\end{align*}
where eigenvalues of $(A-BK_{db})$ are all zero.
The simulation results for the explicit deadbeat MPC closed-loop system is shown in Fig.~\ref{fig: unconstrained terminal eq exp},
which are the same as those in Fig.~\ref{fig: unconstrained terminal eq imp}.

%----------------------------------------------------------
%\subsection{Unconstrained deadbeat MPC with terminal cost}

%In the deadbeat MPC using only terminal cost,
%the implicit solution is calculated by \eqref{unconstrained opt terminal cost} and \eqref{receding horizon scheme},
%and the deadbeat performance is illustrated by Fig.\ref{fig: unconstrained terminal cost imp}.
%It can be seen that performances are the same as those in Figs.~\ref{fig: unconstrained terminal eq imp} and \ref{fig: unconstrained terminal eq exp}.
%
%\begin{figure}
%\begin{center}
%\includegraphics[scale=0.35]{unconstrained_cost_imp.eps}    
%\caption{States and control of the closed-loop system with unconstrained deadbeat MPC using the terminal cost:
%performances are the same as those in Figs.~\ref{fig: unconstrained terminal eq imp} and \ref{fig: unconstrained terminal eq exp}.} 
%\label{fig: unconstrained terminal cost imp}                             
%\end{center}                                
%\end{figure}

%The performances of the explicit deadbeat MPC with only terminal cost are the same with that in Section \ref{sec unconstrained}, %and \ref{sec sim unconstrained eq}, 
%since their control gains are exactly the same.

%----------------------------------------------------------
\subsection{Constrained deadbeat MPC with terminal cost}

Suppose that the control input is subject to hard constraint given by
$-6\leq u\leq 6$.

In the proposed deadbeat MPC, the control horizon is set to $N=3$ which equals the system dimension.
A linear feedback $K =  [2.2150   ~15.0471   ~14.6128]$ is selected, 
such that eigenvalues of $A-BK$ are inside the unit circle.
%Take $Q = I_{3\times3}$ and $R = 0.1$.
The weighting matrix can be solved from the Lyapunov equation \eqref{Lyapunov equation} by
\begin{align*}
	P = \begin{bmatrix}
		 6.1590  & 19.4637   & 5.8132\\
   19.4637 &  96.8173   &40.0964\\
    5.8132  & 40.0964  & 29.9407
	\end{bmatrix}.
\end{align*}

\begin{figure}
\begin{center}
\includegraphics[scale=0.35]{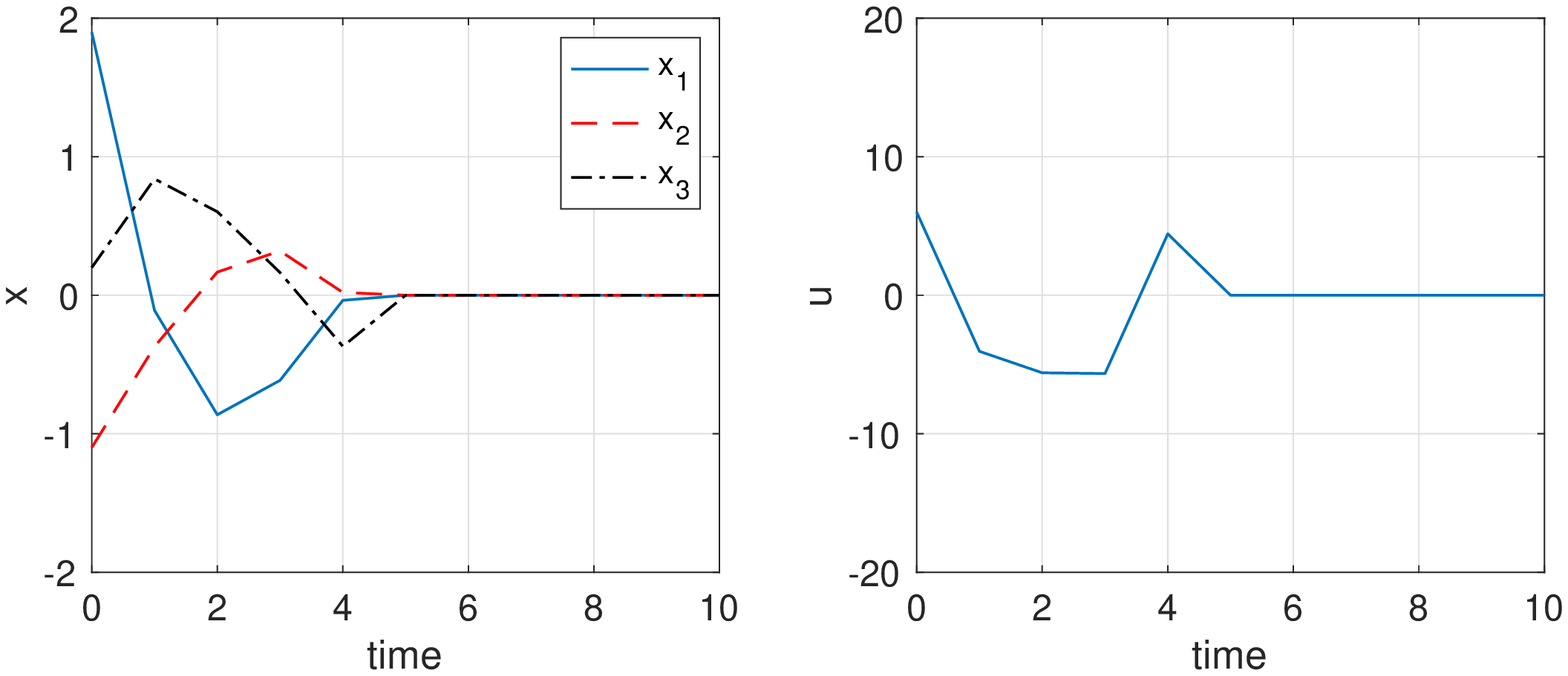}    
\caption{States and control of the closed-loop system with constrained deadbeat MPC:
they converge to the origin in finite time, 
and the transient process is completed in 5 steps.} 
\label{fig: constrained deadbeat}                             
\end{center}                                
\end{figure}

Initial states are assumed to be the same as those in unconstrained cases.
Closed-loop performances are illustrated by Fig.~\ref{fig: constrained deadbeat},
where states and control input is capable of converging to zero in finite time.
The control input is bounded within its constraints,
and it can be seen that the control input during transient process is significantly smaller than that in the unconstrained deadbeat MPC.
Due to the constrained control input, the transient process takes 5 steps, which is relatively longer than that (3 steps) of the unconstrained case.

%====================================================================
\section{Conclusion}

Deadbeat MPC strategies are proposed for SISO linear systems.
The key setting is to assign the control horizon equal to the system dimension.
The proposed unconstrained deadbeat MPC can be designed by either adding an terminal equality constraint or penalizing only the terminal cost.
It is proved that, the unconstrained deadbeat MPC is equialvent to linear deadbeat control.
The proposed constrained deadbeat MPC can be designed by penalizing only the terminal cost subject to the terminal inequality constraint.
It is proved that the closed-loop state is capable of converging to the origin in finite time,
while constraints are always satisfied.

\balance

%-------------------------------------------------------------------------------
%\begin{ack}                               % Place acknowledgements
%%Partially supported by the Roman Senate.  % here.
%This work was supported by National Natural Science Foundation of China under grant no. 62073015.
%\end{ack}

%\bibliographystyle{plain}        % Include this if you use bibtex 

\bibliographystyle{ifacconf}        % Include this if you use bibtex
%\biboptions{authoryear}

\bibliography{deadbeat}           % and a bib file to produce the 
                                 % bibliography (preferred). The
                                 % correct style is generated by
                                 % Elsevier at the time of printing.

%\begin{thebibliography}{99}     % Otherwise use the  
                                 % thebibliography environment.
                                 % Insert the full references here.
                                 % See a recent issue of Automatica 
                                 % for the style.
%  \bibitem[Heritage, 1992]{Heritage:92}
%     (1992) {\it The American Heritage. 
%     Dictionary of the American Language.}
%     Houghton Mifflin Company.
%  \bibitem[Able, 1956]{Abl:56}
%     B.~C.~Able (1956). Nucleic acid content of macroscope. 
%     {\it Nature 2}, 7--9. 
%  \bibitem[Able {\em et al.}, 1954]{AbTaRu:54}   
%     B.~C. Able, R.~A. Tagg, and M.~Rush (1954).
%     Enzyme-catalyzed cellular transanimations.
%     In A.~F.~Round, editor, 
%     {\it Advances in Enzymology Vol. 2} (125--247). 
%     New York, Academic Press.
%  \bibitem[R.~Keohane, 1958]{Keo:58}
%     R.~Keohane (1958).
%     {\it Power and Interdependence: 
%     World Politics in Transition.}
%     Boston, Little, Brown \& Co.
%  \bibitem[Powers, 1985]{Pow:85}
%     T.~Powers (1985).
%     Is there a way out?
%     {\it Harpers, June 1985}, 35--47.

%\end{thebibliography}

\appendix
%\section{A summary of Latin grammar}    % Each appendix must have a short title.
%\section{Some Latin vocabulary}         % Sections and subsections are supported  
                                        % in the appendices.
\end{document}